\documentclass{icm2010}
\usepackage[all]{xy}

\title[Cluster algebras and representation theory]{Cluster algebras and representation theory}
\author[B. Leclerc]{Bernard Leclerc
}
\contact[leclerc@math.unicaen.fr]{LMNO, Universit\'e de Caen, CNRS UMR 6139, F-14032 Caen cedex, France}

\newtheorem{theorem}{Theorem}[section]

\newtheorem{problem}[theorem]{Problem}
\newtheorem{conj}[theorem]{Conjecture}

\theoremstyle{definition}

\def\C{\mathbb C}
\def\CC{\mathcal C}
\def\A{\mathcal A}
\def\Z{\mathbb Z}
\def\N{\mathbb N}
\def\Q{\mathbb Q}
\def\Sl{\mathfrak{sl}}

\def\g{\mathfrak{g}}

\def\n{\mathfrak{n}}
\def\md{\mathrm{mod\,}}

\def\SC{\mathcal S}

\def\F{\mathcal F}
\def\la{\lambda}
\def\ie{\emph{i.e. }}
\def\x{\mathbf x}
\def\y{\mathbf y}
\def\resp{\emph{resp. }}
\def\Si{\Sigma}
\def\a{\alpha}
\def\L{\Lambda}
\def\d{{\mathbf d}}
\def\i{{\mathbf i}}
\def\f{{\mathfrak f}}
\def\M{{\mathcal M}}
\def\de{\delta}
\def\vph{\varphi}
\def\vpi{\varpi}
\def\ZC{{\mathcal Z}}
\def\I{{\mathcal I}}
\def\Ext{\mathrm{Ext}}
\def\End{\mathrm{End}}
\def\Ga{\Gamma}
\def\hP{\widehat{P}}
\def\hla{\widehat{\la}}
\def\hA{\widehat{A}}
\def\hmu{\widehat{\mu}}

\begin{document}

\begin{abstract}
We apply the new theory of cluster algebras of Fomin and Zelevinsky
to study some combinatorial problems arising in Lie theory.
This is joint work with Geiss and Schr\"oer 
(\S\ref{sect3},\,\ref{sect4},\,\ref{sect5},\,\ref{sect6}), 
and with Hernandez (\S \ref{sect8},\,\ref{sect9}).
\end{abstract}



\maketitle


\section{Introduction: two problems in Lie theory}
\label{sect1}


Let $\g$ be a simple complex Lie algebra of type $A$, $D$, or $E$. 
We denote by $G$ a simply-connected complex algebraic group with Lie
algebra $\g$, by $N$ a maximal unipotent subgroup of $G$, by $\n$
its Lie algebra. 
In \cite{Lu1}, Lusztig has introduced the semicanonical basis $\SC$ of the enveloping
algebra $U(\n)$ of $\n$.
Using the duality between $U(\n)$ and the coordinate ring
$\C[N]$ of $N$, one obtains a new basis $\SC^*$ of $\C[N]$ which  
we call the \emph{dual semicanonical basis} \cite{GLS1}.
This basis has remarkable properties. 
For example there is a natural way of
realizing every irreducible finite-dimensional representation
of $\g$ as a subspace $L(\la)$ of $\C[N]$, and 
$\SC^*$ is compatible with this infinite system of subspaces,
that is, $\SC^*\cap L(\la)$
is a basis of $L(\la)$ for every~$\la$.

The definition of the semicanonical basis is geometric (see below \S \ref{sect5}). 
A priori, to describe an element
of $\SC^*$ one needs to compute the Euler characteristics 
of certain complex algebraic varieties.
Here is a simple example in type $A_3$.
Let $V = V_1 \oplus V_2 \oplus V_3$ be a four-dimensional graded vector space
with $V_1 =\C e_1$, $V_2 =\C e_2 \oplus \C e_3$, and $V_3 = \C e_4$. 
There is an element $\vph_X$ of $\SC^*$ attached to the 
nilpotent endomorphism $X$ of $V$ given by  
\[
Xe_1 = e_2,\quad Xe_2 = Xe_3 = 0,\quad Xe_4 = e_3. 
\]
Let $\F_X$ be the variety of complete flags
$F_1 \subset F_2 \subset F_3$ of subspaces of $V$,
which are graded (\emph{i.e.} $F_i = \oplus_j (V_j\cap F_i) \ (1\le i \le 3)$)
and $X$-stable (\emph{i.e.} $XF_i \subset F_i$).
The calculation of $\vph_X$ amounts to computing the Euler
characteristics of the connected components of $\F_X$.
In this case there are four components, two points and two
projective lines, so these Euler numbers are 1,\,1,\,2,\,2.
Unfortunately, such a direct geometric computation
looks rather hopeless in general.
\begin{problem}\label{prob1}
Find a combinatorial algorithm for calculating $\SC^*$. 
\end{problem}

To formulate the second problem we need more notation.
Let $L\g = \g \otimes \C[t,t^{-1}]$ be the loop algebra of $\g$,
and let $U_q(L\g)$ denote the quantum analogue of its enveloping
algebra, introduced by Drinfeld and Jimbo. 
Here we assume that $q\in\C^*$ is not a root of unity.
The finite-dimensional irreducible representations of $U_q(L\g)$
are of special importance because their tensor products give rise to trigonometric $R$-matrices,
that is, to trigonometric solutions of the quantum Yang-Baxter
equation with spectral parameters \cite{J}.  
The question arises whether the tensor product of two given
irreducible representations is again irreducible.
Equivalently, one can ask whether a given irreducible can be 
factored into a tensor product of representations of strictly
smaller dimensions.

For instance, if $\g = \Sl_2$ and $V_n$ is its 
$(n+1)$-dimensional irreducible representation, the loop algebra $L\Sl_2$
acts on $V_n$ by 
\[
 (x\otimes t^k)(v) = z^k xv, \qquad (x\in\Sl_2,\, k\in\Z,\, v\in V_n).
\]
Here $z\in\C^*$ is a fixed number called the evaluation parameter.
Jimbo \cite{J1} has introduced a simple $U_q(L\Sl_2)$-module
$W_{n,z}$, which can be seen as a $q$-analogue of this evaluation
representation. 
Chari and Pressley \cite{CP} have proved that
$W_{n,z} \otimes W_{m,y}$ is an irreducible
$U_q(L\Sl_2)$-module if and only if 
\[
q^{n-m}\frac{z}{y}\not \not\in \left\{q^{\pm(n+m+2-2k)} \mid 0<k\le\min(n,m)\right\}.
\]
In the other direction, they 
showed that every simple object in the category $\md U_q(L\Sl_2)$ of (type 1) finite-dimensional
$U_q(L\Sl_2)$-modules can be written
as a tensor product of modules of the form $W_{n_i,z_i}$ for
some $n_i$ and $z_i$. Thus the modules $W_{n,z}$ can be regarded as the \emph{prime} simple
objects in the tensor category $\md U_q(L\Sl_2)$.

Similarly, for general $\g$ one would like to ask
\begin{problem}\label{prob2}
Find the prime simple objects of $\md U_q(L\g)$, and describe 
the prime tensor factorization of the simple objects. 
\end{problem}

Both problems are quite hard, and we can only offer partial solutions.
An interesting feature is that, in both situations, cluster algebras provide
the natural combinatorial framework to work with.


\section{Cluster algebras}
\label{sect2}


Cluster algebras were invented by Fomin and Zelevinsky \cite{FZ1}
as an abstraction of certain combinatorial structures which they 
had previously discovered while studying total positivity in
semisimple algebraic groups.
A nice introduction \cite{F} to these ideas is given in these proceedings,
with many references to the growing literature on the subject.

A cluster algebra is a commutative ring with a distinguished set
of generators and a particular type of relations. 
Although there can be infinitely many generators and relations,
they are all obtained from a finite number of them
by means of an inductive procedure called \emph{mutation}.

Let us recall the definition.\footnote{For simplicity
we only consider a particular subclass of cluster algebras: 
the antisymmetric cluster algebras of geometric type.
This is sufficient for our purpose.}
We start with the field of rational functions $\F=\Q(x_1,\ldots,x_n)$.
A \emph{seed} in $\F$ is a pair $\Si=(\y,Q)$, where
$\y=(y_1,\ldots,y_n)$ is
a free generating set of $\F$, and 
$Q$ is a quiver (\ie an oriented graph) with vertices
labelled by $\{1,\ldots,n\}$.
We assume that $Q$ has neither loops nor 2-cycles.
For $k=1,\ldots,n$, one defines a new seed 
$\mu_k(\Si)$
as follows. First $\mu_k(y_i) = y_i$ for $i\not = k$,
and
\begin{equation}
\mu_k(y_k)= \frac{\prod_{i\to k}y_i + \prod_{k\to j} y_j}{y_k}, 
\end{equation}
where the first (\resp second) product is over all arrows of 
$Q$ with target (\resp source) $k$.
Next $\mu_k(Q)$ is obtained from $Q$ by
\begin{itemize}
 \item[(a)] adding a new arrow $i\to j$ for every existing pair of arrows
$i\to k$ and $k\to j$;
 \item[(b)] reversing the orientation of every arrow with target or source
equal to $k$;
 \item[(c)] erasing every pair of opposite arrows possibly created by (a).
\end{itemize}
It is easy to check that $\mu_k(\Si)$ is a seed, and $\mu_k(\mu_k(\Si)) = \Si$.
The \emph{mutation class} $\CC(\Si)$ is the set of all seeds obtained 
from $\Si$ by a finite sequence of mutations $\mu_k$. One can think of
the elements of $\CC(\Si)$ as the vertices of an $n$-regular
tree in which every edge stands for a mutation.
If $\Si' = ((y'_1,\ldots,y'_n),Q')$ is a seed in $\CC(\Si)$, 
then the subset $\{y'_1,\ldots,y'_n\}$
is called a \emph{cluster}, and its elements are called \emph{cluster variables}. 
Now, Fomin and Zelevinsky define the \emph{cluster algebra} $\A_\Si$ as the
subring of $\F$ generated by all cluster variables.
Some important elements of $\A_\Si$ are the \emph{cluster monomials},
\ie monomials in the cluster variables supported on a single cluster.

For instance, if $n=2$ and $\Si=((x_1,x_2),Q)$, where $Q$ is the quiver with $a$ 
arrows from $1$ to $2$,
then $\A_\Si$ is the subring of $\Q(x_1,x_2)$ generated 
by the rational functions $x_k$ defined recursively by 
\begin{equation}
x_{k+1}x_{k-1} = 1+ x_k^a,\qquad (k\in \Z).
\end{equation}
The clusters of $\A_\Si$ are the subsets $\{x_k,x_{k+1}\}$,
and the cluster monomials are the special elements of the form
\[
 x_k^l x_{k+1}^m,\qquad (k\in\Z,\ l,m\in \N).
\]
It turns out that when $a=1$, there are only five different clusters
and cluster variables, namely
\[
 x_{5k+1}=x_1,\ \
 x_{5k+2}=x_2,\ \
 x_{5k+3}=\frac{1+x_2}{x_1},\ \
 x_{5k+4}=\frac{1+x_1+x_2}{x_1x_2},\ \
 x_{5k}=\frac{1+x_1}{x_2}.
\]
For $a\ge 2$ though, the sequence $(x_k)$ is no longer periodic and $\A_\Si$
has infinitely many cluster variables. 
 
The first deep results of this theory shown by Fomin and Zelevinsky are:
\begin{theorem}[\cite{FZ1},\cite{FZ2}]\label{thFZ}
\begin{itemize}
\item[\rm (i)] Every cluster variable of $\A_\Si$ is a Laurent polynomial 
with coefficients in $\Z$
in the cluster variables of any single fixed cluster.
\item[\rm (ii)] $\A_\Si$ has finitely many clusters if and only if the mutation
class $\CC(\Si)$ contains a seed whose quiver is an orientation of
a Dynkin diagram of type $A$, $D$, $E$.
\end{itemize}
\end{theorem}

One important open problem \cite{FZ1} is to prove that the coefficients
of the Laurent polynomials in (i) are always positive.
In \S\ref{sect9} below, we give a (conjectural) representation-theoretical
explanation of this positivity for a certain class of cluster algebras. 
More positivity results, based on combinatorial or geometric 
descriptions of these coefficients, have been obtained by 
Musiker, Schiffler and Williams \cite{MSW}, and by Nakajima \cite{N3}.


\section{The cluster structure of $\C{[}N{]}$} 
\label{sect3}


To attack Problem~\ref{prob1} we adopt the following strategy.
We endow $\C[N]$ with the structure
of a cluster algebra\footnote{Here we mean that $\C[N] = \C\otimes_\Z\A$ 
for some cluster algebra $\A$ contained in $\C[N]$.}.
Then we show that all cluster monomials belong to $\SC^*$,
and therefore we obtain a large family of elements of $\SC^*$
which can be calculated by the combinatorial algorithm of 
mutation.

In \cite[\S 2.6]{BFZ} explicit initial seeds for a 
cluster algebra structure in the coordinate ring of the big
cell of the base affine space $G/N$ were described. A simple modification
yields initial seeds for $\C[N]$ (see \cite{GLS3}). 

For instance, if $G=SL_4$ and $N$ is the subgroup of upper unitriangular
matrices, one of these seeds is
\[
((D_{1,2},\,D_{1,3},\,D_{12,23},\,D_{1,4},\,D_{12,34},\,D_{123,234}), Q), 
\]
where $Q$ is the triangular quiver:
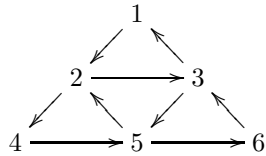
\[
\xymatrix@-1.0pc{
&& \ar[dl]{1} && \\
&\ar[dl]{2}\ar[rr]&&{3}\ar[ul]\ar[dl]\\
{4}\ar[rr]&&{5}\ar[ul]\ar[rr]&&{6}\ar[ul]
} 
\]
Here, by $D_{I,J}$ we mean the regular function on $N$ which
associates to a matrix its minor with row-set $I$
and column-set $J$.
Moreover, the variables 
\[
x_4=D_{1,4},\ x_5=D_{12,34},\ x_6=D_{123,234}
\]
are \emph{frozen}, \ie they cannot be mutated, and therefore they belong
to every cluster. 
Using Theorem~\ref{thFZ}, it is easy to prove that this cluster
algebra has finitely many clusters, namely $14$ clusters and
$12$ cluster variables if we count the $3$ frozen ones.

In general however, that is, for groups $G$ other than $SL_n$ with
$n\le 5$, the cluster structure of $\C[N]$ has infinitely many
cluster variables.
To relate the cluster monomials to $\SC^*$ we have to bring the
preprojective algebra into the picture.


\section{The preprojective algebra}
\label{sect4}


Let $\overline{Q}$ denote the quiver obtained 
from the Dynkin diagram of $\g$ by replacing
every edge by a pair $(\a,\a^*)$
of opposite arrows.
Consider the element
\[
 \rho = \sum (\a\a^* - \a^*\a)
\]
of the path algebra $\C\overline{Q}$ of $\overline{Q}$,
where the sum is over all pairs of opposite arrows.
Following \cite{GP,Ri}, we define the \emph{preprojective 
algebra} $\L$ as the quotient of $\C\overline{Q}$ by 
the two-sided ideal generated by $\rho$.
This is a finite-dimensional selfinjective algebra,
with infinitely many isomorphism classes of indecomposable
modules, except if $\g$ has type $A_n$ with $n\le 4$.
It is remarkable that these few exceptional cases 
coincide precisely with the cases when $\C[N]$ has
finitely many cluster variables. Moreover, it is a nice
exercise to verify that the number of indecomposable 
$\L$-modules is then equal to the number of cluster
variables. 

This suggests a close relationship in general between
$\L$ and $\C[N]$. 
To describe it we start with Lusztig's Lagrangian
construction of the enveloping algebra $U(\n)$ \cite{Lu0, Lu1}.
This is a realization of $U(\n)$ as an algebra of $\C$-valued
constructible functions over the varieties of representations of $\L$.

To be more precise, we need to introduce more notation.
Let $S_i \ (1\le i\le n)$ be the one-dimensional $\L$-modules
attached to the vertices $i$ of $\overline{Q}$. 
Given a sequence $\i=(i_1,\ldots,i_d)$ and a $\L$-module $X$
of dimension $d$,
we introduce the variety $\F_{X,\i}$ of
flags of submodules
\[
\f = (0=F_0\subset F_1 \subset \cdots \subset F_d = X) 
\]
such that $F_k/F_{k-1} \cong S_{i_k}$ for $k=1,\ldots,d$.
This is a projective variety.
Denote by $\L_\d$ the variety of $\L$-modules $X$ with a given
dimension vector $\d=(d_i)$, where $\sum_i d_i =d$.
Consider the constructible function $\chi_\i$ on $\L_\d$  
given by 
\[
\chi_\i(X) = \chi(\F_{X,\i})
\]
where $\chi$ denotes the Euler-Poincar\'e characteristic.
Let $\M_{\d}$ be the $\C$-vector space spanned by
the functions $\chi_\i$ for all possible sequences $\i$
of length $d$, and let 
\[
\M = \bigoplus_{\d\in\N^n} \M_\d. 
\]
Lusztig has endowed $\M$ with an associative multiplication
which formally resembles a convolution product, and he has
shown that, if we denote by $e_i$ the Chevalley generators of $\n$,
there is an algebra isomorphism $U(\n) \stackrel{\sim}{\rightarrow} \M$
mapping the product $e_{i_1}\cdots e_{i_d}$ to $\chi_\i$ for every $\i=(i_1,\ldots,i_d)$.

Now, following \cite{GLS1,GLS2}, we dualize the picture. 
Every $X\in\md\L$ determines a linear form $\de_X$ on $\M$ given by
\[
 \de_X(f) = f(X),\qquad (f\in\M).
\]
Using the isomorphisms $\M^* \simeq U(\n)^* \simeq \C[N]$, the 
form $\de_X$ corresponds to an element $\vph_X$ of $\C[N]$,
and we have thus attached to every object $X$ in $\md\L$
a polynomial function $\vph_X$ on $N$.

For example, if $\g$ is of type $A_3$, and if we denote by $P_i$
the projective cover of $S_i$ in $\md\L$, one has
\[
 \vph_{P_1} = D_{123,234},\quad
 \vph_{P_2} = D_{12,34},\quad
 \vph_{P_3} = D_{1,4}.
\]
More generally, the functions $\vph_X$ corresponding to the
12 indecomposable $\L$-modules are 
the 12 cluster variables of $\C[N]$.

Via the correspondence $X \mapsto \vph_X$ the ring $\C[N]$
can be regarded as a kind of Hall algebra of the category 
$\md\L$. Indeed the multiplication of $\C[N]$ encodes
extensions in $\md\L$, as shown by the following crucial
result. Before stating it, we recall that $\md\L$ possesses
a remarkable symmetry with respect to extensions, namely,
$\Ext^1_\L(X,Y)$ is isomorphic to the dual of 
$\Ext^1_\L(Y,X)$ functorially in $X$ and~$Y$ (see \cite{CB,GLS4}).
In particular $\dim\Ext^1_\L(X,Y) = \dim \Ext^1_\L(Y,X)$ for
every~$X, Y$.

\begin{theorem}[\cite{GLS1,GLS4}]\label{cluster_char}
Let $X, Y \in \md\L$. 
\begin{itemize}
 \item[\rm(i)] We have $\vph_X \vph_Y = \vph_{X\oplus Y}$.
 \item[\rm(ii)] Assume that $\dim \Ext^1_\L(X,Y)=1$, and let
\[
 0\to X \to L \to Y \to 0 
 \quad\mbox{and}\quad 
 0\to Y \to M \to X \to 0
\]
be non-split short exact sequences. Then
$\vph_X \vph_Y = \vph_L + \vph_M$. 
\end{itemize}
\end{theorem}
In fact \cite{GLS4} contains a formula for 
$\vph_X \vph_Y$ valid for any dimension of $\Ext^1_\L(X,Y)$,
but we will not need it here.
As a simple example of (ii) in type $A_2$, one can take 
$X=S_1$ and $Y=S_2$. Then we have the non-split short
exact sequences
\[
 0\to S_1 \to P_2 \to S_2 \to 0 
 \quad\mbox{and}\quad 
 0\to S_2 \to P_1 \to S_1 \to 0,
\]
which imply the relation 
$\vph_{S_1} \vph_{S_2} = \vph_{P_2} + \vph_{P_1}$,
that is, the elementary determinantal relation 
$D_{1,2}D_{2,3} = D_{1,3} + D_{12,23}$
on the unitriangular subgroup of $SL_3$.
More generally, the short Pl\"ucker relations in $SL_{n+1}$
can be obtained as instances of (ii). 

We note that Theorem~\ref{cluster_char} is the analogue for
$\md\L$ of a formula of Caldero and Keller \cite{CK} for 
the cluster categories
introduced by Buan, Marsh, Reineke, Reiten and Todorov \cite{BMRRT}
to model cluster algebras with an acyclic seed.
Cluster categories are not abelian, but Keller \cite{Kel1}
has shown that they are triangulated, so in this setting
exact sequences are replaced by distinguished triangles.


\section{The dual semicanonical basis $\SC^*$}
\label{sect5}


We can now introduce the basis $\SC^*$ of the vector space $\C[N]$.
Let $\d = (d_i)$ be a dimension vector. The variety $\overline{E}_\d$
of representations of $\C\overline{Q}$ with dimension vector $\d$ is a
vector space of dimension $2\sum d_id_j$, where the sum is over all
pairs $\{i,j\}$ of vertices of the Dynkin diagram which are
joined by an edge. This vector space has a natural symplectic structure.
Lusztig \cite{Lu0} has shown that $\L_\d$ is a Lagrangian
subvariety of $\overline{E}_\d$,
and that the number of its irreducible components is equal to 
the dimension of the degree $\d$ homogeneous component of $U(\n)$
(for the standard $\N^n$-grading given by the Chevalley generators). 
Let $\ZC$ be an irreducible component of $\L_\d$. Since the map 
$\vph : X \mapsto \vph_X$
is a constructible map on $\L_\d$, it is constant on a Zariski open subset of $\ZC$.
Let $\vph_\ZC$ denote this generic value of $\vph$ on $\ZC$.
Then, if we denote by $\I = \sqcup_\d \I_\d$ the collection of all
irreducible components of all varieties $\L_\d$, one can easily check
that
\[
 \SC^* = \{\vph_\ZC \mid \ZC \in \I\}
\]
is dual to the basis $\SC = \{f_\ZC \mid \ZC \in \I\}$ of $\M \cong U(\n)$
constructed by Lusztig in \cite{Lu1}, and called by him the semicanonical basis.

For example, if $\g$ is of type $A_n$ and $N$ is the 
unitriangular subgroup in $SL_{n+1}$, then all the
matrix minors $D_{I,J}$ which do not vanish identically on $N$
belong to $\SC^*$ \cite{GLS1}. They are of the form
$\vph_X$, where $X$ is a subquotient of an 
indecomposable projective $\L$-module.

More generally, suppose that $X$ is a \emph{rigid} $\L$-module, 
\ie that $\Ext^1_\L(X,X) = 0$.
Then $X$ is a generic point of the unique irreducible 
component $\ZC$ on which it sits, that is, $\vph_X = \vph_\ZC$
belongs to $\SC^*$, so the calculation of $\vph_\ZC$
amounts to evaluating the Euler characteristics $\chi(\F_{X,\i})$
for every $\i$ (of course only finitely many varieties 
$\F_{X,\i}$ are non-empty).
Thus in type $A_3$, the nilpotent endomorphism $X$ of \S \ref{sect1}
can be regarded as a rigid $\L$-module with dimension vector
$\d = (1,2,1)$, and the connected components of $\F_X$ are
just the non-trivial varieties $\F_{X,\i}$, namely
\[ 
\F_{X,(2,1,2,3)},\quad \F_{X,(2,3,2,1)},\quad \F_{X,(2,2,1,3)},\quad \F_{X,(2,2,3,1)}.
\]
Note however that if $\g$ is not of type $A_n\ (n\le 4)$, there
exist irreducible components $\ZC\in\I$ whose generic points are not 
rigid $\L$-modules.


\section{Rigid $\L$-modules}
\label{sect6}


Let $r$ be the number of positive roots of $\g$.
Equivalently $r$ is the dimension of the affine space $N$.
This is also the number of elements of every cluster of $\C[N]$
(if we include the frozen variables).
Geiss and Schr\"oer have shown \cite{GS} that 
the number of pairwise non-isomorphic indecomposable direct summands
of a rigid $\L$-module is bounded above by $r$.  
A rigid module with $r$ non-isomorphic indecomposable summands 
is called \emph{maximal}.
We will now see that the seeds of the cluster structure of $\C[N]$
come from maximal rigid $\L$-modules.

Let $T=T_1\oplus \cdots \oplus T_r$ be a maximal rigid module,
where every $T_i$ is indecomposable. Define $B = \End_\L T$, a basic
finite-dimensional algebra with simple modules $s_i\ (1\le i\le r)$.
Denote by $\Ga_T$ the quiver of $B$, that is, the quiver with
vertex set $\{1,\ldots,r\}$ and $d_{ij}$ arrows
from $i$ to $j$, where $d_{ij}=\dim \Ext^1_B(s_i,s_j)$.

\begin{theorem}[\cite{GLS2}]
The quiver $\Ga_T$ has no loops nor 2-cycles. 
\end{theorem}

Define $\Si(T) = ((\vph_{T_1},\ldots,\vph_{T_r}),\,\Ga_T)$.

\begin{theorem}[\cite{GLS3}]
There exists an explicit maximal rigid $\L$-module $U$
such that $\Si(U)$ is one of the seeds of the cluster structure
of $\C[N]$. 
\end{theorem}

Let us now lift the notion of seed mutation to the category $\md\L$.

\begin{theorem}[\cite{GLS2}]
Let $T_k$ be a non-projective indecomposable summand of $T$. 
There exists a unique indecomposable module $T_k^*\not\cong T_k$
such that $(T/T_k)\oplus T_k^*$ is maximal rigid. 
\end{theorem}

We call $(T/T_k)\oplus T_k^*$ the \emph{mutation of $T$ in direction $k$},
and denote it by $\mu_k(T)$.
The proof of the next theorem relies among other things on
Theorem~\ref{cluster_char}. 

\begin{theorem}[\cite{GLS2}]\label{thGLS}
\begin{itemize}
\item[{\rm(i)}] We have $\Si(\mu_k(T)) = \mu_k(\Si(T))$, where in the right-hand side
$\mu_k$ stands for the Fomin-Zelevinsky seed mutation.
\item[{\rm(ii)}] The map $T \mapsto \Si(T)$ gives a one-to-one correspondence
between the maximal rigid modules in the mutation class of $U$ and the
clusters of $\C[N]$.
\end{itemize}
\end{theorem}
It follows immediately that the cluster monomials of $\C[N]$ belong to $\SC^*$.
Indeed, by (ii) every cluster monomial is of the form
\[
 \vph_{T_1}^{a_1}\cdots\vph_{T_r}^{a_r} = \vph_{T_1^{a_1}\oplus \cdots \oplus T_r^{a_r}},
\quad (a_1,\ldots,a_r\in\N),
\]
for some maximal rigid module $T=T_1\oplus\cdots\oplus T_r$,
and therefore belongs to $\SC^*$ because $T_1^{a_1}\oplus \cdots \oplus T_r^{a_r}$ is rigid.

Thus the cluster monomials form a large subset of $\SC^*$ which can 
(in principle) be calculated algorithmically by iterating the seed mutation
algorithm from an explicit initial seed. 
This is our partial answer to Problem~\ref{prob1}.

Of course, these results also give a better understanding of the cluster structure
of $\C[N]$. For instance they show immediately that the cluster monomials
are linearly independent (a general conjecture of Fomin and Zelevinsky).
Furthermore, they suggest the definition of new cluster algebra structures
on the coordinate rings of unipotent radicals of parabolic subgroups of $G$,
obtained in a similar manner from some appropriate Frobenius
subcategories of $\md\L$ (see \cite{GLS5}). One can also develop an analogous theory for  
finite-dimensional unipotent subgroups $N(w)$ of a Kac-Moody group attached to
elements $w$ of its Weyl group (see \cite{BIRS,GLS6}).


\section{Finite-dimensional representations of $U_q(L\g)$}
\label{sect7}


We now turn to Problem \ref{prob2}.
We need to recall some known facts about the category 
$\md U_q(L\g)$\footnote{We only consider modules \emph{of type 1}, a mild technical
condition, see \emph{e.g.} \cite[\S12.2 B]{CPbook}.} of finite-dimensional modules over $U_q(L\g)$.

By construction, $U_q(L\g)$ contains a copy of $U_q(\g)$, so in a sense
the representation theory of $U_q(L\g)$ is a refinement of that of
$U_q(\g)$.
Let $\vpi_i\ (1\le i \le n)$ be the fundamental weights of $\g$, and 
denote by
\[
P = \bigoplus_{i=1}^n \Z \vpi_i,\qquad P_+ = \bigoplus_{i=1}^n \N \vpi_i,
\]
the weight lattice and the monoid of dominant integral weights.
It is well known that $\md U_q(\g)$ is a semisimple tensor category,
with simple objects $L(\la)$ parame\-tri\-zed by $\la\in P_+$.
In fact, every $M\in\md U_q(\g)$ has a decomposition
\begin{equation}\label{decompM}
 M = \bigoplus_{\mu \in P} M_\mu
\end{equation}
into eigenspaces for a commutative subalgebra $A$ of $U_q(\g)$ coming
from a Cartan subalgebra of $\g$.
One shows that if $M$ is irreducible, the highest weight occuring
in (\ref{decompM}) is a dominant weight $\la$, $\dim M_\la = 1$, and 
there is a unique simple $U_q(\g)$-module with these properties, hence the notation $M=L(\la)$.
For an arbitrary $M\in \md U_q(\g)$, the formal sum
\[
 \chi(M) = \sum_{\mu \in P} \dim M_\mu\, e^\mu
\]
is called the \emph{character} of $M$, since it characterizes $M$ up to isomorphism.

When dealing with representations of $U_q(L\g)$ one needs to introduce
spectral parameters $z \in \C^*$, and therefore $P$ and $P_+$
have to be replaced by
\[
 \hP = \bigoplus_{1\le i\le n,\ z\in\C^*} \Z (\vpi_i, z),\qquad 
\hP_+ = \bigoplus_{1\le i \le n,\ z\in\C^*} \N (\vpi_i, z).
\]
It was shown by Chari and Pressley \cite{CP,CP2} that finite-dimensional 
irreducible representations
of $U_q(L\g)$ were similarly determined by their highest $l$-weight
$\hla \in \hP_+$
(where $l$ stands for ``loop''). 
This comes from the existence of 
a large commutative subalgebra $\hA$ of $U_q(L\g)$ containing $A$.
If $M\in\md U_q(L\g)$ is regarded as a $U_q(\g)$-module by restriction 
and decomposed as in (\ref{decompM}), 
then every $U_q(\g)$-weight-space
$M_\mu$ has a finer decomposition into generalized eigenspaces for $\hA$
\[
 M_\mu = \bigoplus_{\hmu\in\hP} M_{\hmu}
\]
where the $\hmu = \sum_k m_{i_k}(\vpi_{i_k},z_k)$ in the right-hand side all satisfy
$\sum_k m_{i_k}\vpi_{i_k} = \mu$.
The corresponding formal sum
\[
 \chi_q(M) = \sum_{\hmu \in \hP} \dim M_{\hmu}\, e^{\hmu}
\]
has been introduced by Frenkel and Reshetikhin \cite{FR} and called 
by them the \emph{$q$-character} of $M$. 
It characterizes the class of $M$ in the Grothendieck ring of $\md U_q(L\g)$,
but one should be warned that this is not a semisimple category, so this
is much coarser than an isomorphism class.

For instance, the $4$-dimensional irreducible 
representation $V_3$ of $U_q(\Sl_2)$ with highest weight $\la=3\vpi_1$ has character 
\[
\chi(V_3) = Y^3 + Y^{1} + Y^{-1} + Y^{-3} 
\]
if we set $Y=e^{\vpi_1}$. There is a family
$W_{3,z}\in\md U_q(L\Sl_2)$ of affine analogues  
of~$V_3$, parametrized by $z\in\C^*$, whose $q$-character
is given by
\[
 \chi_q(W_{3,z})= Y_zY_{zq^2}Y_{zq^4}
+ Y_zY_{zq^2}Y_{zq^6}^{-1}
+ Y_zY_{zq^4}^{-1}Y_{zq^6}^{-1}
+ Y_{zq^2}^{-1}Y_{zq^4}^{-1}Y_{zq^6}^{-1},
\]
where we write $Y_a = e^{(\vpi_1,a)}$ for $a\in\C^*$. 
Thus $W_{3,z}$ has highest $l$-weight 
\[
\hla = (\vpi_1,z) + (\vpi_1,zq^2) + (\vpi_1,zq^4).
\] 
The reader can easily imagine what is the general expression
of $\chi_q(W_{n,z})$ for any $(n,z)\in\N\times\C^*$.
It follows that there is a closed formula for the $q$-character of
every finite-dimensional irreducible $U_q(L\Sl_2)$-module since, as already
mentioned, every such module factorizes as a tensor product
of $W_{n_i,z_i}$ and the factors are given by a simple combinatorial rule \cite{CP}.

The situation is far more complicated in general. In particular
it is not always possible to endow an irreducible $U_q(\g)$-module 
with the structure of a $U_q(L\g)$-module. The only general description
of $q$-characters of simple $U_q(L\g)$-modules, due to 
Ginzburg and Vasserot for type $A$ \cite{GV} and to Nakajima in general \cite{N1},
uses intersection cohomology of certain moduli spaces of 
representations of graded preprojective algebras, called graded quiver 
varieties. This yields a Kazhdan-Lusztig type algorithm for
calculating the irreducible $q$-characters \cite{N2}, but this type
of combinatorics does not easily reveal the possible factorizations
of the $q$-characters.


\section{The subcategories $\CC_\ell$}
\label{sect8}


It can be shown that Problem~\ref{prob2} for $\md U_q(L\g)$
can be reduced to the same problem for some much smaller tensor subcategories
$\CC_\ell\ (\ell\in\N)$ which we shall now introduce.

Denote by $L(\hla)$ the simple object of $\md U_q(L\g)$ with
highest $l$-weight $\hla\in\hP_+$.
Since the Dynkin diagram of $\g$ is a tree, it is a bipartite graph.
We denote by $I=I_0\sqcup I_1$ the corresponding partition of the
set of vertices, and we write $\xi_i = 0$ (\resp $\xi_i = 1$) if
$i\in I_0$ (\resp $i\in I_1$).
For $\ell\in\N$, let 
\[
\hP_{+,\ell} = \bigoplus_{1\le i \le n,\ 0\le k \le\ell} \N (\vpi_i, q^{\xi_i + 2k}).
\]
We then define $\CC_\ell$ as the full subcategory of $\md U_q(L\g)$
whose objects $M$ have all their composition factors of the
form $L(\hla)$ with $\hla\in \hP_{+,\ell}$.
It is not difficult to prove \cite{HL} that $\CC_\ell$ is a tensor subcategory,
and that its Grothendieck ring $K_0(\CC_\ell)$ is the polynomial ring in
the $n(\ell + 1)$ classes of fundamental modules
\[
[L(\vpi_i, q^{\xi_i + 2k})], \qquad (1\le i \le n,\ 0\le k \le\ell). 
\]

For example, let $W_{j,a}^{(i)}$ denote the simple object of $\md U_q(L\g)$
with highest $l$-weight
\[
 (\vpi_i,a) + (\vpi_i,aq^2) + \cdots + (\vpi_i, aq^{2j-2}),
\quad (i\in I,\ j\in\N^*,\ a\in\C^*),
\]
a so-called \emph{Kirillov-Reshetikhin module}.
The $q$-characters of the Kirillov-Resheti\-khin modules satisfy 
a nice system of recurrence relations, called $T$-system in the
physics literature, which allows to calculate them inductively
in terms of the $q$-characters of the fundamental modules
$L(\vpi_i,a)$. 
This was conjectured by Kuniba, Nakanishi and Suzuki \cite{KNS}, 
and proved by Nakajima \cite{N4} (see also \cite{H2} for the non simply-laced
cases).
The $q$-characters of the fundamental modules can in turn be calculated
by means of the Frenkel-Mukhin algorithm \cite{FM}.
One should therefore regard the Kirillov-Reshetikhin modules as the most
``accessible'' simple $U_q(L\g)$-modules.
There are $n(\ell+1)(\ell+2)/2$ such modules in $\CC_\ell$, namely:
\[
W^{(i)}_{j,q^{\xi_i + 2k}},
\quad
(i\in I,\ 0<j\le \ell+1,\ 0\le k \le \ell +1 - j).
\]
 

\section{The cluster algebras $\A_\ell$}
\label{sect9}


Let $Q$ denote the quiver obtained by orienting the Dynkin
diagram of $\g$ so that every $i\in I_0$ (\resp $i\in I_1$)
is a source (\resp a sink).
We define a new quiver $\Ga_\ell$ with vertex set 
$\{(i,k)\mid i\in I,\ 1\le k \le \ell + 1\}$.
There are three types of arrows
\begin{itemize}
 \item[(a)] arrows $(i,k) \rightarrow (j,k)$ for every arrow $i \to j$ in $Q$ 
and every $1\le k\le \ell + 1$;
 \item[(b)] arrows $(j,k) \rightarrow (i,k+1)$ for every arrow $i \to j$ in $Q$ 
and every $1\le k\le \ell$;
 \item[(c)] arrows $(i,k) \leftarrow (i,k+1)$ for every $i\in I$ and every 
$1\le k\le \ell$.
\end{itemize}
For example, if $\g$ has type $A_3$ and $I_0 = \{1,3\}$, the quiver $\Ga_3$ is:

\medskip
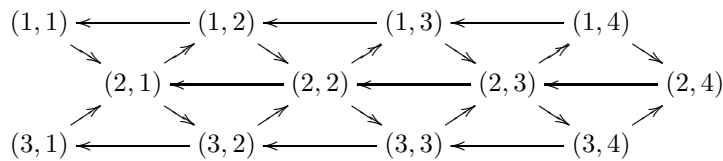
\[
\xymatrix@-1.4pc{
(1,1)\ar[rd] && \ar[ll](1,2) \ar[rd] && \ar[ll](1,3) \ar[rd] && \ar[ll]{(1,4)} \ar[rd] \\
&(2,1)\ar[ru]\ar[rd]&&\ar[ll](2,2)\ar[ru]\ar[rd]&&\ar[ll](2,3)\ar[ru]\ar[rd]&&\ar[ll]{(2,4)}& \\
(3,1)\ar[ru] && \ar[ll](3,2) \ar[ru] && \ar[ll](3,3) \ar[ru] && \ar[ll]{(3,4)} \ar[ru]
}
\]
Let $\x = \{x_{(i,k)} \mid i\in I,\ 1\le k \le \ell + 1\}$ be a set of indeterminates
corresponding to the vertices of $\Ga_\ell$,
and consider the seed $(\x , \Ga_\ell)$ in which the $n$ variables
$x_{(i,\ell+1)}\ (i\in I)$ are frozen.
This is the initial seed of a cluster algebra $\A_\ell \subset \Q(\x)$.
By Theorem~\ref{thFZ}, if $\g$ has type $A_1$ then
$\A_\ell$ has finite cluster type $A_\ell$. Also, if $\ell = 1$, 
$\A_\ell$ has finite cluster type equal to the Dynkin type of $\g$.
Otherwise, except for a few small rank cases, $\A_\ell$ has infinitely
many cluster variables.

Our partial conjectural solution of Problem~\ref{prob2} can be summarized
as follows (see \cite{HL} for more details):
\begin{conj}\label{conjec}
There is a ring isomorphism 
$\iota_\ell : \A_{\ell} \stackrel{\sim}{\to} K_0(\CC_\ell)$
such that
\[
\iota_\ell(x_{(i,k)}) = \left[W^{(i)}_{k,\,q^{{\xi_i}+2({\ell}+1-k)}}\right],\qquad (i\in I,\ 1\le k\le \ell+1). 
\]
The images by $\iota_\ell$ of the cluster variables are classes
of prime simple modules, and the images of the cluster monomials
are the classes of all \emph{real} simple modules in~$\CC_\ell$, \ie those simple
modules whose tensor square is simple. 
\end{conj}
Thus, if true, Conjecture~\ref{conjec} gives a combinatorial description in terms of
cluster algebras of the prime tensor factorization of every real
simple module.
Note that, by definition, the square of a cluster monomial is again a cluster monomial.
This explains why cluster monomials can only correspond to real simple modules.
For $\g=\Sl_2$, all simple $U_q(L\g)$-modules are real.
However for $\g\not = \Sl_2$ there exist \emph{imaginary} simple $U_q(L\g)$-modules
(\ie simple modules whose tensor square is not simple), as shown
in \cite{L}.
This is consistent with the expectation that 
a cluster algebra with infinitely many cluster variables
is not spanned by its set of cluster monomials.

We arrived at Conjecture~\ref{conjec} by noting that the $T$-system
equations satisfied by Kirillov-Reshetikhin modules are 
of the same form as the cluster exchange relations.
This was inspired by the seminal work of Fomin and Zelevinsky \cite{FZ0},
in which cluster algebra combinatorics is used to prove Zamolodchikov's 
periodicity conjecture for $Y$-systems attached to Dynkin diagrams.
Kedem \cite{Ked} and Di Francesco \cite{KDF}, Keller \cite{Kel2,Kel3},
Inoue, Iyama, Kuniba, Nakanishi and Suzuki \cite{IIKNS}, have
also exploited the similarity between cluster exchange relations
and other types of functional equations arising in mathematical physics
($Q$-systems, generalized $T$-systems,
$Y$-systems attached to pairs of simply-laced Dynkin diagrams).
Recently, Inoue, Iyama, Keller, Kuniba and Nakanishi \cite{IIKKN1,IIKKN2}
have obtained a proof of the periodicity conjecture for all $T$-systems and $Y$-systems
attached to a non simply-laced quantum affine algebra.

As evidence for Conjecture~\ref{conjec}, we can easily check that for $\g=\Sl_2$
and any $\ell\in\N$, it follows from the results 
of Chari and Pressley \cite{CP}.
On the other hand, for arbitrary $\g$ we have:
\begin{theorem}[\cite{HL,N3}]\label{th9.2}
Conjecture~\ref{conjec} holds for $\g$ of type $A, D, E$ and $\ell = 1$. 
\end{theorem}
This was first proved in \cite{HL} for type $A$ and $D_4$ by combinatorial
and represen\-ta\-tion-theoretic methods, and soon after, by Nakajima \cite{N3} 
in the general case, by using the geometric description of the simple $U_q(L\g)$-modules.
In both approaches, a crucial part of the proof can be summarized in
the following chart:
\[
\begin{matrix}
\mbox{$F$-polynomials} & \leftrightarrow & \mbox{quiver Grassmannians} \cr
\updownarrow && \updownarrow \cr
\mbox{$q$-characters} & \leftrightarrow & \mbox{Nakajima quiver varieties}
\end{matrix}
\]
Here, the $F$-polynomials are certain polynomials 
introduced by Fomin and Zelevinsky \cite{FZ4} which allow
to calculate the cluster variables in terms of a fixed
initial seed.
By work of Caldero-Chapoton \cite{CC}, Fu-Keller \cite{FK} and Derksen-Weyman-Zelevinsky
\cite{DWZ1,DWZ2}, $F$-polynomials have a geometric description via Grassmannians
of subrepresentations of some quiver representations attached to
cluster variables: this is the
upper horizontal arrow of our diagram.
The lower horizontal arrow refers to the already mentioned 
relation between irreducible $q$-characters and
perverse sheaves on quiver varieties established by Nakajima 
\cite{N1,N2}. In \cite{HL} we have shown that the $F$-polynomials
for $\A_1$ are equal to certain natural truncations of the 
corresponding irreducible $q$-characters of $\CC_1$ (the left vertical arrow), 
and we observed that this yielded an alternative geometric description of
these $q$-characters in terms of ordinary homology of quiver Grassmannians.
In \cite{N3} Nakajima used a Deligne-Fourier transform to obtain a direct 
relation between perverse sheaves on quiver varieties for $\CC_1$
and homology of quiver Grassmannians (the right vertical arrow), and deduced from it
the desired connection with the cluster algebra $\A_1$.

The other main step in the approach of \cite{HL} is a certain tensor product
theorem for the category $\CC_1$. It states that a tensor product
$S_1\otimes \cdots \otimes S_k$ of simples objects of $\CC_1$ is simple
if and only if $S_i\otimes S_j$ is simple for every pair $1\le i < j\le k$.
A generalization of this theorem to the whole category $\md U_q(L\g)$
has been recently proved by Hernandez \cite{H}. Note that the theorem of Hernandez is
also valid for non simply-laced Lie algebras $\g$, and thus opens the way to a similar
treatment of Problem~\ref{prob2} in this case.

Conjecture~\ref{conjec} has also been checked 
for $\g$ of type $A_2$ and $\ell = 2$ \cite[\S 13]{HL}. 
In that small rank case, $\A_2$ still has finite cluster type $D_4$, 
and this implies that $\CC_2$ has only real objects.
There are 18 explicit prime simple objects with respective dimensions
\[
3,\ 3,\ 3,\ 3,\ 3,\ 3,\ 6,\ 6,\ 6,\ 6,\ 8,\ 8,\ 8,\ 10,\ 10,\ 15,\ 15,\ 35,  
\]
and 50 factorization patterns (corresponding to the 50 vertices of
a generalized associahedron of type $D_4$ \cite{FZ2}).
Our proof in this case is quite indirect and uses a lot of
ingredients: the quantum affine Schur-Weyl duality,
Ariki's theorem for type $A$ affine Hecke algebras \cite{A}, the
coincidence of Lusztig's dual canonical and dual semicanonical bases
of $\C[N]$ in type $A_4$ \cite{GLS1}, and Theorem~\ref{thGLS}.

One remarkable consequence of Theorem~\ref{th9.2}
from the point of view of cluster algebras is that it immediately implies
the positivity conjecture of Fomin and Zelevinsky for the cluster algebras
$\A_1$ with respect to any reference cluster (see \cite[\S2]{HL}).
Conjecture~\ref{conjec} would similarly yield positivity for the
whole class of cluster algebras $\A_\ell$.


\section{An intriguing relation} 
\label{sect10}


Problem~\ref{prob1} and Problem~\ref{prob2} may not be 
as unrelated as it would first seem.
For a suggestive example, let us take $\g$ of type $A_3$. 
In that case, the abelian category $\md \L$ has 12 indecomposable 
objects (which are all rigid), 3 of them being projective-injective.
On the other hand the tensor category $\CC_1$ has 12 prime simple objects
(which are all real), 3 of them having the property that their tensor product
with every simple of $\CC_1$ is simple.
It is easy to check that $\C[N]$ and $\C\otimes_\Z K_0(\CC_1)$ are isomorphic
as (complexified) cluster algebras with frozen variables.
Therefore we have a unique one-to-one correspondence 
\[
 X\ \leftrightarrow\ S
\]
between rigid 
objects $X$ of $\md \L$ and simple objects $S$ of $\CC_1$ 
such that
\[
 \vph_X \equiv [S],
\]
that is, such that $X$ and $S$ project to the same cluster monomial.
In this correspondence, direct sums $X \oplus X'$ map to 
tensor products $S\otimes S'$.
It would be interesting to find a general framework for relating 
in a similar way, via cluster algebras, certain additive categories 
such as $\md \L$ to certain tensor categories such as $\CC_1$.
We refer to \cite{K} for a very accessible survey of these ideas.
  
\bigskip\noindent
\emph{Acknowledgements}: I wish to thank Christof Geiss and Jan Schr\"oer 
for their very precious collaboration and friendship over the years,
and David Hernandez for his stimulating enthusiasm. 


\end{document}